\documentclass[10pt]{article}

\usepackage{gitinfo}

\makeatletter
\def\cahierline{\gdef\@cahierline}
\cahierline{\textcolor{gray}{Cahier du GERAD G-\number\year-\cahiernumber}}
\def\gitline{\gdef\@gitline}
\gitline{\textcolor{gray}{Commit \gitAbbrevHash\ by \gitAuthorName\ on \gitAuthorIsoDate}}
\def\ps@mytitlepage{\let\@mkboth\@gobbletwo
  \def\@oddfoot{\normalfont\scriptsize\ttfamily\hfil\@gitline}
  \def\@evenfoot{\normalfont\scriptsize\ttfamily\@gitline\hfil}
  \def\@evenhead{\normalfont\small\rlap{\@cahierline}%\hspace{\headlineindent}
  \hfil}
  \def\@oddhead{\normalfont\small\hfil%\hspace{\headlineindent}
  \llap{\@cahierline}}
}

\def\ps@myheadings{\let\@mkboth\@gobbletwo
  \def\@oddfoot{\normalfont\scriptsize\ttfamily\@cahierline\hfil\@gitline}
  \def\@evenfoot{\normalfont\scriptsize\ttfamily\@gitline\hfil\@cahierline}
  \def\@evenhead{\normalfont\small\rlap{\thepage}%\hspace{\headlineindent}%
  }
  \def\@oddhead{\normalfont\small\scshape%\hyperlink{contents}{[toc]}
  \hfil\rightmark%\hspace{\headlineindent}
  \llap{\thepage}}
}
\makeatother

\newcommand{\cahiernumber}{12}  % Insert your Cahier du GERAD number.

\usepackage[margin=1in]{geometry}
\usepackage[T1]{fontenc}
\usepackage{amsthm}  % must come before newpxtext and newpxmath
\usepackage{amssymb}

\usepackage{eulervm} % Euler math

\usepackage{amsmath}
\usepackage{thmtools}
\usepackage[final]{microtype} % Less badboxes
\usepackage{xspace}
\usepackage{subdepth}  % Equalize depth of subscripts. Amazing!

\usepackage{graphicx}
\graphicspath{{figs/}}

\usepackage{xcolor}
\usepackage{tikz}
\usepackage{tikzscale}
\usetikzlibrary{external}
\newcommand*{\includetikzgraphics}[2][]{%
  \tikzsetnextfilename{#2}%
  \includegraphics[#1]{#2.tikz}
}

\usepackage{pgfplots}
\pgfplotsset{compat=newest}
\usepackage{algorithmicx}
\usepackage[section]{algorithm}
\usepackage[noend]{algpseudocode}

\usepackage[theme=grayscale,charsperline=90]{jlcode}

\usepackage[pdfpagelabels]{hyperref}   % Internal hyperlinks
\hypersetup{
  breaklinks,
  colorlinks,
  linkcolor={blue!70!black},
  urlcolor={red!50!black},
  citecolor={green!50!black},
  pdfborder={0 0 0},      % No borders around internal hyperlinks
  pdfauthor={Dominique Orban} % authors
}
\usepackage[nameinlink]{cleveref}

\declaretheorem[name=Lemma]{lemma}
\declaretheorem[name=Proposition]{proposition}
\declaretheorem[name=Corollary]{corollary}

\makeatletter
\def\@maketitle{%
  \newpage
  \vspace*{-\topskip}      % remove the initial space
  \begingroup\centering    % instead of \begin{center}
  \let\footnote\thanks
  \hrule height \z@        % to avoid the insertion of lineskip glue
    {\LARGE \@title \par}%
    \vskip 1.5em
    {\large
      \lineskip .5em
      \begin{tabular}[t]{c}%
        \@author
      \end{tabular}\par}%
    \vskip 1em
    \@date%
  \par\endgroup            % instead of \end{center}
  \vskip 1.5em             % <--- modify this to adjust the separation
}
\makeatother

\usepackage[numbers,sort&compress]{natbib}
\makeatletter
\renewcommand\NAT@bibsetnum[1]{\settowidth\labelwidth{\@biblabel{#1}}%
   \setlength{\leftmargin}{\bibindent}\addtolength{\leftmargin}{\dimexpr\labelwidth+\labelsep\relax}%
   \setlength{\itemindent}{-\bibindent}%
   \setlength{\listparindent}{\itemindent}
\setlength{\itemsep}{\bibsep}\setlength{\parsep}{\z@}%
   \ifNAT@openbib
     \addtolength{\leftmargin}{\bibindent}%
     \setlength{\itemindent}{-\bibindent}%
     \setlength{\listparindent}{\itemindent}%
     \setlength{\parsep}{0pt}%
   \fi
}
\makeatother

\usepackage{etoolbox}
\makeatletter
\patchcmd{\NAT@test}{\else \NAT@nm}{\else \NAT@hyper@{\NAT@nm}}{}{}
\makeatother
\usepackage{dsfont}
\newcommand{\R}{\mathds{R}}
\newcommand{\N}{\mathds{N}}

\newcommand{\B}{\mathds{B}}

\newcommand{\minim}{\mathop{\textup{minimize}}}
\newcommand{\minimize}[1]{\minim_{#1}}
\newcommand{\st}{\textup{subject to}}
\newcommand{\argminim}{\mathop{\textup{argmin}}}
\newcommand{\argmin}[1]{\argminim_{#1}}
\newcommand{\argmaxim}{\mathop{\textup{argmax}}}
\newcommand{\argmax}[1]{\argmaxim_{#1}}
\newcommand{\supp}{\mathop{\textup{supp}}}
\newcommand{\proj}{\mathop{\textup{proj}}}

\usepackage{newunicodechar}
\newunicodechar{Δ}{\ensuremath{\Delta}}

\title{Computing a Sparse Projection into a Box}
\author{%
  Dominique Orban\footnote{%
    GERAD and Department of Mathematics and Industrial Engineering, Polytechnique Montr\'eal. E-mail: \href{mailto:dominique.orban@gerad.ca}{dominique.orban@gerad.ca}.
  }
  \thanks{Research supported by an NSERC Discovery grant.}%
}

\begin{document}
  \maketitle
  \thispagestyle{mytitlepage}

  \begin{abstract}
    We describe a procedure to compute a projection of \(w \in \R^n\) into the intersection of the so-called \emph{zero-norm} ball \(k \B_0\) of radius \(k\), i.e., the set of \(k\)-sparse vectors, with a box centered at a point of \(k \B_0\).
    The need for such projection arises in the context of certain trust-region methods for nonsmooth regularized optimization.
    Although the set into which we wish to project is nonconvex, we show that a solution may be found in \(O(n \log(n))\) operations.
    We describe our Julia implementation and illustrate our procedure in the context of two trust-region methods for nonsmooth regularized optimization.
  \end{abstract}

  % \begin{resume}
  %   Nous proposons une procédure pour calculer une projection de \(w \in \R^n\) dans l'intersection de la soi-disant boule en \emph{norme zéro} \(k \B_0\) de rayon \(k\), c'est-à-dire l'ensemble des vecteurs ayant au plus \(k\) composantes non nulles, et d'une boîte centrée en point de \(k \B_0\).
  %   Cette projection est nécessaire dans le contexte de certaines méthodes de région de confiance pour l'optimisation non lisse régularisée.
  %   Bien que l'ensemble dans lequel on projette est non convexe, il est possible d'obtenir une solution en \(O(n \log(n))\) opérations.
  %   Nous décrivons notre implémentation dans le langage Julia et illustrons la procédure dans le context de deux méthodes de région de confiance pour l'optimisation non lisse régularisée.
  % \end{resume}

  \pagestyle{myheadings}

  \section{Introduction}

  We describe a procedure to compute a projection of a point in \(\R^n\) into the intersection of the set of \(k\)-sparse vectors with a box centered at a \(k\)-sparse vector.

  Specifically, let \(\Delta \B_\infty\) be the \(\ell_{\infty}\)-norm ball of radius \(\Delta \geq 0\) and centered at the origin, and \(x + \Delta \B_\infty\) be the same ball centered at \(x \in \R^n\).
  The set of \(k\)-sparse vectors in \(\R^n\), otherwise known as the \(\ell_0\)-pseudonorm ``ball'' of radius \(k \in \{0, 1, \ldots, n\}\), is denoted \(k \B_0\) and is the set of vectors with at most \(k\) nonzero components.
  Assume that \(x \in k \B_0\).
  For given \(w \in \R^n\), we seek to compute
  \begin{equation}%
    \label{eq:proj}
    p(w) \in P(w) := \argminim \ \{ \|w - y\|_2 \mid y \in C \}
    \qquad
    C := k \B_0 \cap (x + \Delta \B_\infty).
  \end{equation}
  Because \(C\) is closed, \(P(w) \neq \varnothing\), but because \(C\) is nonconvex, \(P(w)\) may contain several elements.
  In~\eqref{eq:proj}, we seek a global minimum---local nonglobal minima sometimes exist, but are of no particular interest here.
  Although it may appear as though the problem has exponential complexity due to the combinatorial nature of \(k\)-sparsity, we show that a solution may be found in \(O(n \log(n))\) operations.
  We describe our Julia implementation and illustrate our procedure in the context of two trust-region methods for nonsmooth regularized optimization.

  \subsection*{Context}

  The computation of~\eqref{eq:proj} occurs in the evaluation of proximal operators encountered during the iterations of the trust-region method of \citet{aravkin-baraldi-orban-2021} for nonsmooth regularized optimization.
  Their method is designed for problems of the form
  \begin{equation}%
    \label{eq:minf+h}
    \minimize{x \in \R^n} \ f(x) + h(x),
  \end{equation}
  where \(f: \R^n \to \R\) has Lipschitz-continuous gradient and \(h: \R^n \to \R \cup \{ \pm \infty \}\) is lower semi-continuous and proper.
  In large-scale data fitting and signal reconstruction problems, \(h(x) = \chi(x \mid k \B_0)\) encodes sparsity constraints and is of interest if one is to recover a solution with at most \(k\) nonzero elements, where \(\chi(\cdot \mid A)\) is the indicator of \(A \subseteq \R^n\), i.e.,
  \[
    \chi(x \mid A) =
    \begin{cases}
      0 & \text{if } x \in A, \\
      \infty & \text{otherwise}.
    \end{cases}
  \]
  All iterates \(x_j\) generated are feasible in the sense that \(x_j \in k \B_0\).
  At iteration \(j\), a step \(s\) is computed in
  \[
    \argmin{u} \tfrac{1}{2} \|u - v\|_2^2 + h(x_j + u) + \chi(u \mid \Delta_j \B_\infty),
  \]
  where \(v \in \R^n\) is given and \(\Delta_j \B_\infty\) is the trust region centered at the origin of radius \(\Delta_j > 0\).
  With the change of variables \(z := x_j + u\), we may rewrite the above as
  \[
    \argmin{z} \ \tfrac{1}{2} \|z - w\|_2^2 + \chi(z \mid k \B_0) + \chi(z \mid x_j + \Delta_j B_\infty) - \{x_j\},
  \]
  where \(w := x_j + v\), which precisely amounts to~\eqref{eq:proj} with \(x_j\) in the role of \(x\) and \(\Delta_j\) in the role of \(\Delta\) because the two indicators may be combined into the indicator of the intersection.

  Because nonsmooth regularized problems often involve a nonlinear least squares smooth term, \citet{aravkin-baraldi-orban-2021b} develop a Levenberg-Marquardt variant of their trust region method.
  The latter requires the same projections as just described.

  \subsection*{Notation}

  Let \(\supp(x) := \{i = 1, \ldots, n \mid x_i \neq 0\}\) be the \emph{support} of \(x\).
  If \(A \subseteq \R^n\) is closed and \(A \neq \varnothing\), we denote
  \[
    \proj(w \mid A) := \argminim \ \{ \|w - y\|_2 \mid y \in A \},
  \]
  the projection of \(w\) into \(A\), which is a set with at least one element.

  When the projection of \(w\) into \(A\) is unique, such as happens when \(A\) is convex, we slightly abuse notation and write
  \(y = \proj(w \mid A)\) instead of \(\{y\} = \proj(w \mid A)\).

  If \(B \subseteq \R^n\), the notation \(\proj(\proj(w \mid A) \mid B)\) refers to the set \(\{z \in \proj(y \mid B) \text{ for some } y \in \proj(w \mid A)\}\).

  If \(S \subseteq \{1, \ldots, n\}\), the cardinality of \(S\) is denoted \(|S|\), and its complement is \(S^c\).
  For such \(S\) and for \(x \in \R^n\), we denote \(x_S\) the subvector of \(x\) indexed by \(S\) and \(A_S := \{x \in \R^n \mid x_{S^c} = 0\}\).
  Clearly, \(0 \in A_S\) for any such \(S\).

  Because
  \[
    k \B_0 = \bigcup \ \{ A_S \mid S \subseteq \{1, \ldots, n\}, \ |S| = k \},
  \]
  \citep[p.~\(175\)]{beck-2017}, we refer to \(A_S\) as a \emph{piece} of \(k \B_0\).

  \subsection*{Related Research}

  \citet{duchi-2008} describe how to project efficiently into the \(\ell_1\)-norm ball.
  The \(\ell_1\)-norm is probably the most widely used convex approximation of the \(\ell_0\) norm as minimizing \(\|x\|_1\) promotes sparsity under certain conditions---see, e.g., \citep{candes-romberg-tao-2006} and the vast ensuing compressed sensing literature.
  % By contrast, minimizing \(\|x\|_0\) is NP hard in general.

  \citet{gupta2010l1} describe how to project into the intersection of an \(\ell_1\)-norm ball with a box, which may be seen as a relaxation of~\eqref{eq:proj}.
  \citet{thom-palm-2013} and \citet{thom-rapp-palm-2015} propose a linear-time and constant space algorithm to compute a projection into a hypersphere with a prescribed sparsity, where sparsity is measured by the ratio of the \(\ell_1\) to the \(\ell_2\) norm.

  \citet{beck-eldar-2013} provide optimality conditions for the minimization of a smooth function over \(k \B_0\).
  \citet{beck-hallak-2016} provide optimality conditions for problems of the form~\eqref{eq:proj} where the box is replaced with a symmetric set satisfying certain conditions.
  Unfortunately,~\eqref{eq:proj} does not satisfy those conditions unless \(x = 0\), at which point it is easy to see that a solution simply consists in chaining the projection into \(k \B_0\) with that into \(\Delta \B_\infty\).
  That is what \citet[Proposition~\(4.3\)]{luss-teboulle-2013} do with \(1\B_2\) instead of \(\Delta \B_\infty\).

  \citet[Proposition~\(4\)]{bolte-sabach-teboulle-2014} show how to project into the intersection of \(k \B_0\) with the nonnegative orthant.

  \citet{pmlr-v28-kyrillidis13} explain how to compute a sparse projection into the simplex, which is probably the most closely related research to our objectives.
  The simplex necessarily intersects all pieces of \(k \B_0\), which need not be the case for~\eqref{eq:proj}.

  \section{Geometric Intuition}

\begin{figure}[ht]
    \centering
    \begin{tikzpicture}[scale=0.75]
      % Axis
      \draw[thin] (-4, 0) -- (4, 0);
      \draw[thin] (0, -5) -- (0, 3);

      \draw[green, ultra thick] (-3, 0) -- (3, 0);
      \draw[green, ultra thick] (0, 2) -- (0, -4);

      \fill[black] (0, -1) circle (2pt) node[right] {$x$};
      \draw[very thick] (-3, 2) -- (3, 2) -- (3, -4) -- (-3, -4) -- (-3, 2);

      \fill[black] (-2.5, 2.5) circle (2pt) node[above right] {$w_1$};
      \draw (0, 0) -- (-3, 3);
      \draw[red, ->, very thick] (-2.5, 2.5) -- (0, 2.5);
      \draw[red, ->, very thick] (0, 2.5) -- (0, 2);
      \draw[red, ->, very thick] (-2.5, 2.5) -- (-2.5, 0);

      \fill[black] (-4, -3.5) circle (2pt) node[below] {$w_2$};
      \draw[red, ->, very thick] (-4, -3.5) -- (-3, -3.5);
      \draw[red, ->, very thick] (-3, -3.5) -- (0, -3.5);
      \draw[blue, ->, very thick] (-4, -3.5) -- (-4, 0);
      \draw[blue, ->, very thick] (-4, 0) -- (-3, 0);

      \fill[black] (-3, 0) circle (2pt) node[above left]   {$p(w_2)$};
      \fill[black] (-2.5, 0) circle (2pt) node[below right] {$p(w_1)$};
      \fill[black] (0, 2) circle (2pt) node[above right] {$p_2$};
    \end{tikzpicture}
    \qquad
    \begin{tikzpicture}[scale=0.75]
      % Axis
      \draw[thin] (-4, 0) -- (4, 0);
      \draw[thin] (0, -5) -- (0, 3);

      \draw[green, ultra thick] (0, -1) -- (0, -4);

      \fill[black] (0, -2.5) circle (2pt) node[right] {$x$};
      \draw[very thick] (-1.5, -1) -- (1.5, -1) -- (1.5, -4) -- (-1.5, -4) -- (-1.5, -1);

      \fill[black] (3, 1) circle (2pt) node[right] {$w$};
      \fill[black] (0, -1) circle (2pt) node[above left] {$p_1 \in p(w)$};

      \draw[red, ->, very thick] (3, 1) -- (3, 0);
      \draw[red, ->, very thick] (3, 0) -- (1.5, -1);
      \draw[blue, ->, very thick] (3, 1) -- (1.5, -1);
      \draw[blue, ->, very thick] (1.5, -1) -- (1.5, 0);
    \end{tikzpicture}
    \caption{%
      Simply composing the projection into \(1 \B_0\) with that into \(x + \Delta \B_\infty\), in either order, may lead to an erroneous projection.
      \label{fig:2projs}
    }
  \end{figure}
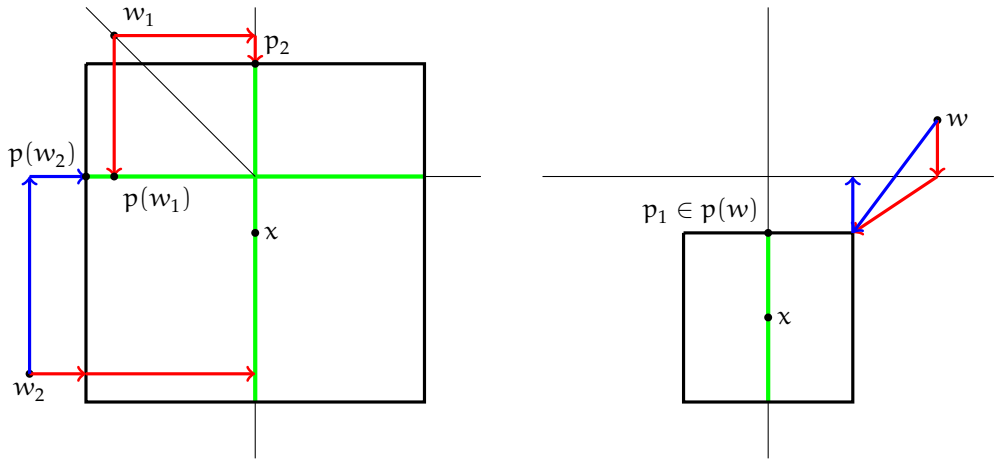

  Naively chaining the projection into one set with that into the other, in either order, does not necessarily yield a point into the intersection of the two sets, even if the latter are convex.
  \Cref{fig:1B0,fig:2projs} illustrates two situations that we may encounter when \(k = 1\) and \(n = 2\).

  A few simple observations about \Cref{fig:1B0,fig:2projs} reveal some difficulties associated with the computation of \(p(w)\):
  \begin{enumerate}
    \item because both components of \(w_1\) are equal in absolute value, as indicated by the thin diagonal in \Cref{fig:2projs}, \(\proj(w_1 \mid k \B_0)\) is a set with two elements, and projecting those into \(x + \Delta \B_\infty\) yields \(p(w_1)\) (the correct global minimum) and \(p_2\) (a spurious local minimum);
    \item moving \(w_1\) up slightly would preserve \(p(w_1)\), but projecting into \(1 \B_0\) first would lead to \(p_2\);
    \item moving \(w_1\) slightly to the right would result in a projection that is slightly to the right of \(p(w_1)\) on the figure, but projecting into \(1 \B_0\) first would lead to \(p_2\);
    \item moving \(w_1\) further to the right would result in \(P(w_1) = \{p(w_1), p_2\}\) and moving it further still would result in \(P(w_1) = \{p_2\}\);
    \item in the rightmost plot, chaining the projections either way leads to a point that does not even lie in the intersection.
  \end{enumerate}

  \begin{figure}[ht]
    \centering
    \begin{tikzpicture}[scale=0.75]
      % Axis
      \draw[thin] (-4, 0) -- (4, 0);
      \draw[thin] (0, -5) -- (0, 3);

      \draw[green, ultra thick] (-3, 0) -- (3, 0);
      \draw[green, ultra thick] (0, 2) -- (0, -4);

      \fill[black] (0, -1) circle (2pt) node[right] {$x$};
      \draw[very thick] (-3, 2) -- (3, 2) -- (3, -4) -- (-3, -4) -- (-3, 2);

      \fill[black] (-2.5, 2.5) circle (2pt) node[above right] {$w_1$};
      \fill[black] (-4, -3.5) circle (2pt) node[below] {$w_2$};
      \fill[black] (3, 2) circle (2pt) node[right] {$w_3$};

      \fill[black] (-3, 0) circle (2pt) node[above left]   {$p_1$};
      \fill[black] (-2.5, 0) circle (2pt) node[below right] {$p(w_1)$};
      \fill[black] (0, 2) circle (2pt) node[above right] {$p_2$};
      \fill[black] (3, 0) circle (2pt) node[above right] {$p_3 \in P(w_3)$};
      \fill[black] (0, -4) circle (2pt) node[above right] {$p_4$};
    \end{tikzpicture}
    \qquad
    \begin{tikzpicture}[scale=0.75]
      % Axis
      \draw[thin] (-4, 0) -- (4, 0);
      \draw[thin] (0, -5) -- (0, 3);

      \draw[green, ultra thick] (0, -1) -- (0, -4);

      \fill[black] (0, -2.5) circle (2pt) node[right] {$x$};
      \draw[very thick] (-1.5, -1) -- (1.5, -1) -- (1.5, -4) -- (-1.5, -4) -- (-1.5, -1);

      \fill[black] (2, 1) circle (2pt) node[right] {$w$};
      \fill[black] (0, -1) circle (2pt) node[above left] {$p_1 \in p(w)$};
      \fill[black] (0, -4) circle (2pt) node[above right] {$p_2$};
    \end{tikzpicture}
    \caption{%
      \label{fig:1B0}
      The set composed of the two axes is \(1 \B_0\) in \(\R^2\), the box is \(x + \Delta \B_\infty\) and the green set is their intersection.
      Left: \(P(w_1) = \{p(w_1)\}\), and \(P(w_2) = \{p_1\}\).
      With respect to \(w_1\), the other cardinal points are \(p_2\), a local minimum, \(p_3\), a local maximum, and \(p_4\), a global maximum.
      % However, \(P(w_3) = \{p_2, p_3\}\).
      Right: the intersection of \(1 \B_0\) with \(x + \Delta \B_\infty\) is entirely determined by \(\supp(x)\), \(P(w) = \{p_1\}\) while \(p_2\) is a global maximum.
    }
  \end{figure}
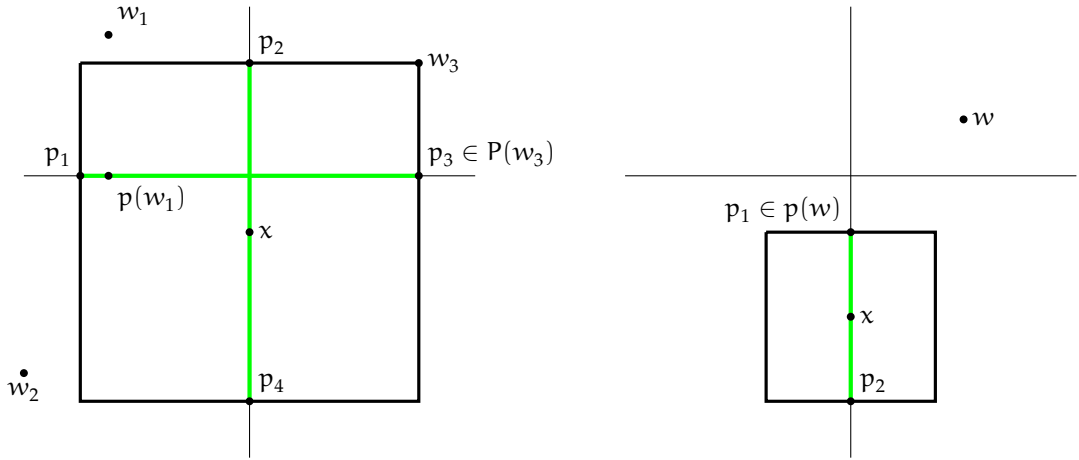

  Note that \(1 \B_0\) is a special case for any value of \(n\): its intersection with \(x + \Delta \B_\infty\) consists in either a single line segment, or \(n\) segments.
  Indeed, the first possibility is that the nonzero component of \(x\) is \(|x_i| > \Delta\).
  In that case, any \(y \in 1 \B_0\) with \(y_j \neq 0\) and \(i \neq j\) satisfies \(\|y - x\|_\infty \geq |x_i| > \Delta\), and therefore \(y \not \in x + \Delta \B_\infty\).
  The only other possibility is that \(|x_i| \leq \Delta\), in which case \(0 \in x + \Delta \B_\infty\), and therefore, all pieces of \(1 \B_0\) intersect the box.

  For \(1 < k < n\), however, the intersection may consist in any number of pieces between \(1\) and \({n \choose k}\).

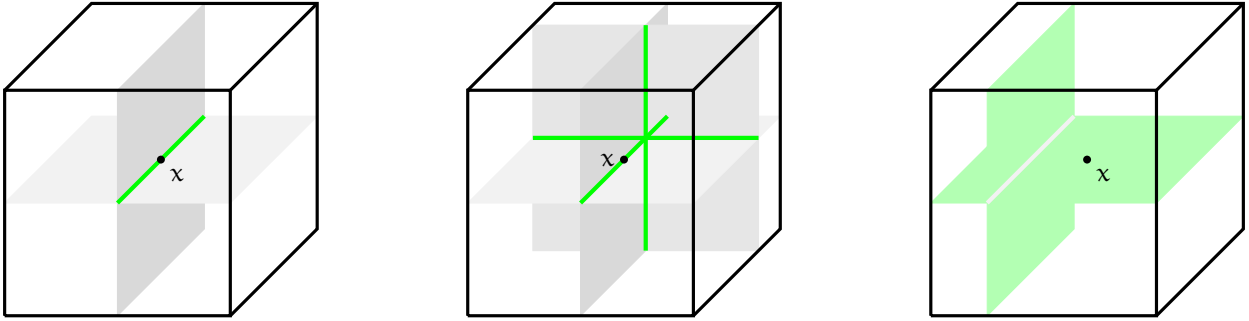
\begin{figure}[ht]
    \centering
    % !!! The TikZ coordinate system is (y, z, x) for some stupid reason
    \begin{tikzpicture}[scale=0.75]
      % draw one half plane at a time in the right order to show which are in the front and which are in the back
      % left part of xy plane
      \draw[gray!10, fill=gray!10] (-2, 0, 5) -- (0, 0, 5) -- (0, 0, 1) -- (-2, 0, 1) -- (-2, 0, 5);
      % xz plane
      \draw[gray!30, fill=gray!30] (0, -2, 5) -- (0, -2, 1) -- (0, 2, 1) -- (0, 2, 5) -- (0, -2, 5);
      % right part of xy plane
      \draw[gray!10, fill=gray!10] (0, 0, 5) -- (2, 0, 5) -- (2, 0, 1) -- (0, 0, 1) -- (0, 0, 5);

      \draw[green, ultra thick] (0, 0, 1) -- (0, 0, 5);
      \fill[black] (0, 0, 3) circle (2pt) node[below right] {$x$};

      % front face
      \draw[very thick] (-2, -2, 5) -- (2, -2, 5) -- (2, 2, 5) -- (-2, 2, 5) -- (-2, -2, 5);
      % back face
      \draw[very thick] (-2, 2, 5) -- (-2, 2, 1) -- (2, 2, 1) -- (2, -2, 1) -- (2, -2, 5);
      % missing edge
      \draw[very thick] (2, 2, 5) -- (2, 2, 1);
    \end{tikzpicture}
    \hfill
    \begin{tikzpicture}[scale=0.75]
      % draw one part of a plane at a time in the right order to show which are in the front and which are in the back
      % bottom left part of yz plane
      \draw[gray!20, fill=gray!20] (0, 0, 0) -- (-2, 0, 0) -- (-2, -2, 0) -- (0, -2, 0) -- (0, 0, 0);
      % left part of xy plane
      \draw[gray!10, fill=gray!10] (-2, 0, 3) -- (0, 0, 3) -- (0, 0, -1) -- (-2, 0, -1) -- (-2, 0, 3);
      % top left part of yz plane
      \draw[gray!20, fill=gray!20] (0, 2, 0) -- (-2, 2, 0) -- (-2, 0, 0) -- (0, 0, 0) -- (0, 2, 0);
      % xz plane
      \draw[gray!30, fill=gray!30] (0, -2, 3) -- (0, -2, -1) -- (0, 2, -1) -- (0, 2, 3) -- (0, -2, 3);
      % bottom right part of yz plane
      \draw[gray!20, fill=gray!20] (0, 0, 0) -- (2, 0, 0) -- (2, -2, 0) -- (0, -2, 0) -- (0, 0, 0);
      % right part of xy plane
      \draw[gray!10, fill=gray!10] (0, 0, 3) -- (2, 0, 3) -- (2, 0, -1) -- (0, 0, -1) -- (0, 0, 3);
      % top right part of yz plane
      \draw[gray!20, fill=gray!20] (0, 0, 0) -- (0, 2, 0) -- (2, 2, 0) -- (2, 0, 0) -- (0, 0, 0);

      \draw[green, ultra thick] (0, 0, -1) -- (0, 0, 3);
      \draw[green, ultra thick] (-2, 0, 0) -- (2, 0, 0);
      \draw[green, ultra thick] (0, -2, 0) -- (0, 2, 0);
      \fill[black] (0, 0, 1) circle (2pt) node[left] {$x$};

      % front face
      \draw[very thick] (-2, -2, 3) -- (2, -2, 3) -- (2, 2, 3) -- (-2, 2, 3) -- (-2, -2, 3);
      % back face
      \draw[very thick] (-2, 2, 3) -- (-2, 2, -1) -- (2, 2, -1) -- (2, -2, -1) -- (2, -2, 3);
      % missing edge
      \draw[very thick] (2, 2, 3) -- (2, 2, -1);
    \end{tikzpicture}
    \hfill
    \begin{tikzpicture}[scale=0.75]
      \draw[green!30, fill=green!30] (-1, 0, 5) -- (3, 0, 5) -- (3, 0, 1) -- (-1, 0, 1) -- (-1, 0, 5);
      \draw[green!30, fill=green!30] (0, -2, 5) -- (0, -2, 1) -- (0, 2, 1) -- (0, 2, 5) -- (0, -2, 5);
      \draw[gray!10, ultra thick] (0, 0, 1) -- (0, 0, 5);

      \fill[black] (1, 0, 3) circle (2pt) node[below right] {$x$};

      % front face
      \draw[very thick] (-1, -2, 5) -- (3, -2, 5) -- (3, 2, 5) -- (-1, 2, 5) -- (-1, -2, 5);
      % back face
      \draw[very thick] (-1, 2, 5) -- (-1, 2, 1) -- (3, 2, 1) -- (3, -2, 1) -- (3, -2, 5);
      % missing edge
      \draw[very thick] (3, 2, 5) -- (3, 2, 1);

    \end{tikzpicture}
    \caption{%
      \label{fig:B0inR3}
      Left: the green segment represents a possible intersection of \(1 \B_0\) with a box in \(\R^3\).
      The gray plane sections only serve to position the segment visually in three dimensions.
      Center: another possible intersection of \(1 \B_0\) with a box in \(\R^3\).
      The box either intersects a single axis, or all of them.
      Right: the green region is a possible intersection of \(2 \B_0\) with a box in \(\R^3\).
      The gray segment only serves as a visual aid and is part of the intersection.
    }
  \end{figure}

  \Cref{fig:B0inR3} illustrates situations that may arise for \(k = 1\) or \(2\) and \(n = 3\).

  \section{Background and Preliminary Results}

  The unique projection \(y\) of any \(w\) into \(x + \Delta \B_\infty\) has components
  \[
    y_i = \max(x_i - \Delta, \min(w_i, x_i + \Delta)), \quad i = 1, \ldots, n.
  \]

  Given \(S \subseteq \{1, \ldots, n\}\), we obtain the unique projection of any \(w\) into \(A_S\) by setting \(w_i = 0\) for all \(i \in S^c\).

  A projection \(y\) of any \(w\) into \(k \B_0\) is a vector that has the same \(k\) largest components in absolute value as \(w\), and the rest of its components set to zero \citep[Lemma~\(6.71\)]{beck-2017}.

  In the vein of \citet{beck-eldar-2013}, it is possible to state necessary optimality conditions for the more general problem
  \begin{equation}%
    \label{eq:min-sparse-box}
    \minimize{y \in \R^n} \ f(y) \quad \st \ y \in C,
  \end{equation}
  of which~\eqref{eq:proj} is a special case.
  Despite the fact that our algorithm is not based on such necessary conditions, they are relevant in their own right, and we now review and specialize them to~\eqref{eq:proj}.

  \begin{lemma}%
    \label{lem:bf}
    Let \(y^\star\) be a solution of~\eqref{eq:min-sparse-box} where \(f\) is continuously differentiable.
    \begin{enumerate}
      \item If \(\|y^\star\|_0 < k\), then for all \(i = 1, \ldots, n\),
      \[
        \frac{\partial f(y^\star)}{\partial y_i}
        \begin{cases}
          \leq 0 & \text{ if } y^\star_i = x_i + \Delta \\
          \geq 0 & \text{ if } y^\star_i = x_i - \Delta \\
          = 0 & \text{ otherwise;}
        \end{cases}
      \]
      \item if \(\|y^\star\|_0 = k\), the same conditions hold for all \(i \in \supp(y^\star)\).
    \end{enumerate}
  \end{lemma}

  \begin{proof}
    The proof follows that of \citep[Theorem~\(2.1\)]{beck-eldar-2013}.
    If \(\|y^\star\|_0 < k\), then for all \(i = 1, \ldots, n\),
    \[
      0 \in \argmin{t \in \R} \ \{g(t) \mid \|y^\star + t e_i - x\|_\infty \leq \Delta\},
    \]
    where \(e_i\) is the \(i\)-th column of the identity, and \(g(t) := f(y^\star + t e_i)\).

    Because \(y^\star \in x + \Delta \B_\infty\), the constraint above reduces to \(|y^\star_i + t - x_i| \leq \Delta\).
    The conclusion follows directly from the standard KKT conditions by noting that \(g'(0) = \partial f(y^\star) / \partial y_i\).

    If \(\|y^\star\|_0 = k\), the same reasoning goes for all \(i \in \supp(y^\star)\).
  \end{proof}

  By analogy with \citep[Theorem~\(2.1\)]{beck-eldar-2013}, a candidate satisfying the conditions of \Cref{lem:bf} is called a \emph{basic feasible} point.
  The following corollary follows directly from \Cref{lem:bf} with \(f(y) := \tfrac{1}{2} \|w - y\|_2^2\).

  \begin{corollary}%
    \label{cor:bf}
    Let \(y^\star\) be a solution of~\eqref{eq:proj}.
    \begin{enumerate}
      \item If \(\|y^\star\|_0 < k\), then for all \(i = 1, \ldots, n\),
      \[
        y^\star_i
        \begin{cases}
          \leq w_i & \text{ if } y^\star_i = x_i + \Delta \\
          \geq w_i & \text{ if } y^\star_i = x_i - \Delta \\
          = w_i & \text{ otherwise;}
        \end{cases}
      \]
      \item if \(\|y^\star\|_0 = k\), the same conditions hold for all \(i \in \supp(y^\star)\).
    \end{enumerate}
  \end{corollary}

  \Cref{lem:bf} and \Cref{cor:bf} are only necessary conditions, and they are rather weak; there often exist vectors satisfying the conditions stated that are not solutions of~\eqref{eq:min-sparse-box} or~\eqref{eq:proj}.
  Consider for example \(k = 1\) in \(\R^2\), \(x = (0, -1)\), \(\Delta = 2\), and \(w = (2, 3)\).
  Then, \(y = (2, 0)\) satisfies the conditions of \Cref{cor:bf}: \(\|y\|_0 = 1\), \(\supp(y) = \{1\}\) and \(y_1 = x_1 + \Delta \leq w_1\).
  However, \(P(w) = \{(0, 1)\}\).
  Indeed, \(\|w - (0, 1)\| = 2 \sqrt{2} < 3 = \|w - y\|\).

  Observe that thanks to \citep[Lemma~\(2.1\)]{beck-eldar-2013}, the number of basic feasible points of~\eqref{eq:proj} is finite.
  Therefore, so is the cardinality of \(P(w)\).

  For a constant \(L > 0\), \citeauthor{beck-eldar-2013} define \(y \in C\) to be \(L\)-stationary for~\eqref{eq:min-sparse-box} if it satisfies \(y \in \proj(y - L^{-1} \nabla f(y) \mid C)\), a condition insipired by optimality conditions for convex problems.
  They state the following result, whose proof remains valid for~\eqref{eq:min-sparse-box}.

  \begin{lemma}[{\protect \citealp[Lemma~\(2.2\)]{beck-eldar-2013}}]%
    \label{lem:L-stat}
    For any \(L > 0\), \(y \in \R^n\) is \(L\)-stationary for~\eqref{eq:min-sparse-box} if and only if \(y \in C\) and
    \[
      \frac{\partial f(y)}{\partial y_i} = 0 \quad (i \in \supp(y))
      \qquad \text{and} \qquad
      \left| \frac{\partial f(y)}{\partial y_i} \right| \leq L M_k(y) \quad (i \not \in \supp(y)),
    \]
    where \(M_k(y)\) is the \(k\)th largest component of \(y\) in absolute value.
  \end{lemma}

  With \(f(y) := \tfrac{1}{2} \|w - y\|_2^2\), \(L\)-stationarity reads \(y \in \proj(y - L^{-1} (y - w) \mid C)\).
  Due to the simple form of \(\nabla f(y) = y - w\), \Cref{lem:L-stat} specializes as follows.

  \begin{corollary}%
    \label{cor:L-stat}
    For any \(L > 0\), \(y \in \R^n\) is \(L\)-stationary for~\eqref{eq:proj} if and only if \(y \in C\) and
    \[
      w_i = y_i \quad (i \in \supp(y)),
      \qquad \text{and} \qquad
      |w_i| \leq L M_k(y) \quad (i \not \in \supp(y)).
    \]
  \end{corollary}

  As a special case of \Cref{cor:L-stat}, if \(\|y\|_0 < k\), then \(M_k(y) = 0\) and we obtain \(w_i = 0\) for \(i \not \in \supp(y)\).
  In that case, \(L\)-stationarity turns out to be independent of \(L\) and requires that \(y = w\), i.e., there is a unique \(L\)-stationary point if \(w \in C\), and there are no \(L\)-stationary points if \(w \not \in C\).

  \(L\)-stationarity is stronger than basic feasibility in the sense that if \(y\) is \(L\)-stationary for~\eqref{eq:min-sparse-box} for any \(L > 0\), then \(y\) is also a basic feasible point \citep[Corollary\(2.1\)]{beck-eldar-2013}.

  Under a Lipschitz assumption, solutions of~\eqref{eq:min-sparse-box} are \(L\)-stationary, as stated in the following result.

  \begin{proposition}[{\protect \citealp[Theorem~\(2.2\)]{beck-eldar-2013}}]%
    \label{prop:L-stat}
    Assume \(\nabla f\) is Lipschitz continuous with constant \(L_f\) and \(y\) solves~\eqref{eq:min-sparse-box}.
    Then, for any \(L > L_f\),
    \begin{enumerate}
      \item \(y\) is \(L\)-stationary;
      \item \(\proj(y - L^{-1} \nabla f(y) \mid C)\) is a singleton.
    \end{enumerate}
  \end{proposition}

  \begin{proof}
    The proof follows by verifying that \citep[Lemma~\(2.4\)]{beck-eldar-2013} continues to hold for~\eqref{eq:min-sparse-box} and the proof of \citep[Theorem~\(2.2\)]{beck-eldar-2013} holds unchanged.
  \end{proof}

  \Cref{prop:L-stat} clearly applies to~\eqref{eq:proj} as the gradient of \(f(y) := \tfrac{1}{2} \|w - y\|_2^2\) is Lipschitz continuous with constant \(L_f = 1\).
  Thus, solutions of~\eqref{eq:proj} are \(L\)-stationary for \(L > 1\).

  Based on \(L\)-stationarity, \citet{beck-eldar-2013} study the iteration \(y^+ \in \proj(y - L^{-1} \nabla f(y) \mid C)\) and show convergence to an \(L\)-stationary point for~\eqref{eq:min-sparse-box} under the assumption that \(\nabla f\) is Lipschitz continuous.
  Unfortunately, in the case of~\eqref{eq:proj}, solving the subproblem is as difficult as solving~\eqref{eq:proj} directly.

  Finally, \citeauthor{beck-eldar-2013} define the concept of componentwise (CW) optimality as follows: \(y \in C\) is CW-minimum for~\eqref{eq:min-sparse-box} if
  \begin{enumerate}
    \item \(\|y\|_0 < k\) and \(f(y) = \min_t f(y + t e_i)\) for \(i = 1, \dots, n\), or
    \item \(\|y\|_0 = k\) and \(f(y) \leq \min_t f(y - y_i e_i + t e_j)\) for \(i \in \supp(y)\) and \(j = 1, \dots, n\).
  \end{enumerate}
  They observe that any solution is a CW-minimum \citep[Theorem~\(2.3\)]{beck-eldar-2013} and that any CW-minimum is a basic feasible point \citep[Lemma~\(2.5\)]{beck-eldar-2013}.
  The concept of CW-minimum allows them to show that any solution of~\eqref{eq:min-sparse-box} is \(L\)-stationary for a value \(L\) that can be significantly smaller than \(L_f\).
  Based on those observations, they propose two coordinate descent-type methods that converge to a CW-minimum.

  In the next section, we present a number of properties of~\eqref{eq:proj} and an algorithm that identifies a solution directly, without resort to the above stationarity conditions.

  \section{Computing the Projection}

  We begin with a few simple observations.

  \begin{lemma}%
    \label{lem:projAS}
    Let \(S \subseteq \{1, \ldots, n\}\) such that \(|S| = k\).
    If \(y \in x + \Delta \B_\infty\) and \(z = \proj(y \mid A_S)\), then \(\|z - x\|_\infty \leq \|y - x\|_\infty\) and, in particular, \(z \in x + \Delta \B_\infty\).
  \end{lemma}

  \begin{proof}
    Without loss of generality, we may write \(z = (y_S, 0)\).
    Observe now that
    \[
      \Delta \geq
      \|y - x\|_\infty =
      \max(\|y_S - x_S\|_\infty, \|y_{S^c} - x_{S^c}\|_{\infty}) =
      \max(\|y_S - x_S\|_\infty, \|y_{S^c}\|_{\infty}) \geq
      \|y_S - x_S\|_\infty =
      \|z - x\|_\infty,
    \]
    because \(x_{S^c} = 0\).
  \end{proof}

  \Cref{lem:projAS} holds because of the geometry of \(k \B_0\) respective to \(\B_\infty\) and is specific to the \(\ell_\infty\)-norm.
  Indeed, consider for example a ball defined in the \(\ell_2\)-norm and set \(x = (0, -1) \in 1 \B_0\) and \(\Delta = 2\).
  For \(y_1 = (\tfrac{3}{4}, -\tfrac{1}{4})\), we have \(z_1 = \proj(y_1 \mid 1 \B_0) = (\tfrac{3}{4}, 0)\) and \(\|z_1 - x\|_2 > \|y - x\|_2\).
  In this example, \(z_1 \in x + \Delta \B_2\), but consider now \(y_2 = (2, -1)\).
  Then, \(z_2 = \proj(y_2 \mid 1 \B_0) = (2, 0) \not \in x + \Delta \B_2\).

  \begin{lemma}%
    \label{lem:proj1proj2}
    If \(w \in x + \Delta \B_\infty\), then \(P(w) = \proj(w \mid k \B_0)\).
  \end{lemma}

  \begin{proof}
    Any \(y \in \proj(w \mid k \B_0)\) has the same \(k\) largest components in absolute value as \(w\), and the rest of its components set to zero.
    Thus, there must exist \(S \subseteq \{1, \ldots, n\}\) with \(|S| = k\) such that \(y = \proj(w \mid A_S)\).
    By \Cref{lem:projAS}, \(\|y - x\|_\infty \leq \|w - x\|_\infty \leq \Delta\) so that \(y \in x + \Delta \B_\infty\), and hence, \(y \in C\).
    If there were \(z \in C\) such that \(\|z - w\|_2 < \|y - w\|_2\), because \(z \in k \B_0\), there would be a contradiction with the definition of \(y\).
    Therefore, \(y\) is a closest point to \(w\) in \(C\).
  \end{proof}

  % In general, if \(w \in k \B_0\), \(p(w) \in P(w)\) need not have any element in common with \(\proj(w \mid x + \Delta \B_\infty)\).
  % Consider for instance \(x = (3, 0) \in 1 \B_0\), \(\Delta = 2\) and \(w = (0, 1) \in 1 \B_0\).
  % Then, \(proj(w \mid x + \Delta \B_\infty) = (1, 1) \not \in 1 \B_0\).

  \begin{lemma}%
    \label{lem:suppx}
    \(C = A_{\supp(x)} \cap (x + \Delta \B_\infty)\) if and only if \ \(|x_i| > \Delta\) for all \(i \in \supp(x)\).
  \end{lemma}

  \begin{proof}
    The result follows from the observation that for any \(i \in \supp(x)\), there is no \(y \in x + \Delta \B_\infty\) with \(y_i = 0\).
    Indeed, if \(y_i = 0\), \(\|y - x\|_\infty \geq |y_i - x_i| = |x_i| > \Delta\).
  \end{proof}

  \begin{lemma}%
    \label{lem:proj-suppx}
    For any \(S \subseteq \{1, \ldots, n\}\) and any \(w \in \R^n\),
    \[
      \proj(w \mid A_S \cap (x + \Delta \B_\infty)) =
      \proj(proj(w \mid x + \Delta \B_\infty) \mid A_S) =
      \proj(proj(w \mid A_S) \mid x + \Delta \B_\infty),
    \]
    whose unique element is the vector \(y\) such that \(y_S = \proj(w_S \mid x_S + \Delta \B_\infty)\) and \(y_{S^c} = 0\).

    In particular, if \(C = A_{\supp(x)} \cap (x + \Delta \B_\infty)\), then \(P(w) = \{\proj(\proj(w \mid x + \Delta \B_\infty) \mid A_{\supp(x)})\}\).
  \end{lemma}

  \begin{proof}
    The projection is unique because \(A_S \cap (x + \Delta \B_\infty)\) is convex.
    If \(y := \proj(w \mid x + \Delta \B_\infty)\) observe that \(z := \proj(y \mid A_S) \in A_S \cap (x + \Delta \B_\infty)\) by \Cref{lem:projAS}.

    In order to show that \(z \in \proj(w \mid A_S \cap (x + \Delta \B_\infty))\), pick any other \(\bar{z} \in A_S \cap (x + \Delta \B_\infty)\).
    By construction, \(\bar{z} = (\bar{y}_S, 0)\) for some \(\bar{y}\).
    Among the infinitely many possible \(\bar{y}\), we may choose the one such that \(\bar{y}_{S^c} = y_{S^c}\).
    Then,
    \[
      \|w - y\|_2^2 =
      \|w_S - y_S\|_2^2 + \|w_{S^c} - y_{S^c}\|_2^2 =
      \|w_S - z_S\|_2^2 + \|w_{S^c} - y_{S^c}\|_2^2,
    \]
    and
    \[
      \|w - \bar{y}\|_2^2 =
      \|w_S - \bar{y}_S\|_2^2 + \|w_{S^c} - \bar{y}_{S^c}\|_2^2 =
      \|w_S - \bar{z}_S\|_2^2 + \|w_{S^c} - y_{S^c}\|_2^2.
    \]
    By definition of \(y\), \(\|w - y\|_2 \leq \|w - \bar{y}\|_2\) and the above therefore implies \(\|w_S - z_S\|_2^2 \leq \|w_S - \bar{z}_S\|_2^2\).
    Because \(z_{S^c} = \bar{z}_{S^c} = 0\), we may add \(\|w_{S^c}\|_2^2\) to both sides of the previous inequality to obtain \(\|w - z\|_2 \leq \|w - \bar{z}\|_2\).
  \end{proof}

  \Cref{lem:suppx} provides an easily computable criterion to determine that \(C = A_{\supp(x)} \cap (x + \Delta \B_\infty)\), and, thanks to \Cref{lem:proj-suppx}, we find an element of \(P(w)\) by setting all components of \(\proj(w \mid x + \Delta \B_\infty)\) that are not in \(\supp(x)\) to zero.
  Such situation is represented in the rightmost plot of \Cref{fig:1B0}.

  By \Cref{lem:suppx}, if there is \(|x_i| \leq \Delta\), then \(x + \Delta \B_\infty\) intersects other pieces of \(k \B_0\) than \(A_{\supp(x)}\).
  We now determine which pieces, and their number.
  Let
  \[
    s(x) := \{ i \in  \supp(x) \mid |x_i| \leq \Delta\}
    \quad \text{and} \quad
    \ell(x) := \{ i \in \supp(x) \mid |x_i| > \Delta \}
  \]
  be the \emph{small} and \emph{large} nonzero components of \(x\).

  In the special case where \(s(x) = \supp(x)\), i.e., \emph{all} nonzero components of \(x\) are small, \(x + \Delta \B_\infty\) intersects \emph{all} pieces of \(k \B_0\) because \(0 \in C\).
  Unfortunately, there are
  \[
    {n \choose k} =
    \frac{n!}{k! \, (n - k)!}
  \]
  of them.
  As it turns out, it is possible to compute \(p(w) \in P(w)\) for any \(w \in \R^n\) in \(O(n \log(n))\) operations.
  In view of \Cref{lem:proj1proj2}, we assume that \(w \not \in x + \Delta \B_\infty\).

  % Let us first return temporarily to the case where \(s(x) = \supp(x)\).
  We may decompose~\eqref{eq:proj} as suggested in \citep{pmlr-v28-kyrillidis13} and observe that \(y_\star \in \proj(w \mid C)\) if and only if \(S_\star\) and \(y_\star\) are in
  \begin{equation}%
    \label{eq:proj-split}
    \argmin{\substack{S \subseteq \{1, \ldots, n\}\\ |S| = k}} \ \argmin{y \in A_S \cap (x + \Delta \B_\infty)} \ \|w - y\|_2^2.
  \end{equation}
  In the case of \(\B_\infty\), we know that \(y \in A_S \cap (x + \Delta \B_\infty)\) if and only if \(y \in x + \Delta \B_\infty\) and \(y_{S^c} = 0\), i.e., if and only if \(y_S \in x_S + \Delta \B_\infty\) and \(y_{S^c} = 0\).
  Thus, we may rewrite~\eqref{eq:proj-split} as
  \[
    \argmin{\substack{S \subseteq \{1, \ldots, n\}\\ |S| = k}} \ \argmin{\substack{y \in x + \Delta \B_\infty\\ y_{S^c} = 0}} \ \|w_S - y_S\|_2^2 + \|w_{S^c}\|^2 =
    \argmin{\substack{S \subseteq \{1, \ldots, n\}\\ |S| = k}} \ \argmin{\substack{y_S \in x_S + \Delta \B_\infty\\ y_{S^c} = 0}} \ \|w_S - y_S\|_2^2 - \|w_S\|^2.
  \]
  For fixed \(S\), the unique solution of the inner problem is \(y = y(S)\) such that \(y_S = \proj(w_S \mid x_S + \Delta \B_\infty)\) and \(y_{S^c} = 0\).
  Thus, the problem reduces to finding the optimal piece, determined by
  \begin{equation}%
    \label{eq:proj-S}
    S_\star \in \argmax{\substack{S \subseteq \{1, \ldots, n\}\\ |S| = k}} \ \|w_S\|^2 - \|w_S - y_S\|^2.
  \end{equation}
  Because~\eqref{eq:proj-S} requires examining all pieces of \(k \B_0\), it may be solved by noting that
  \[
    \|w_S\|^2 - \|w_S - y_S\|^2 = e^T z, \quad
    e = (1, 1, \ldots, 1), \quad
    z_i = w_i^2 - {(w_i - y_i)}^2, \ i \in S,
  \]
  i.e., the objective is the sum of the components of \(z\) with indices in \(S\).
  Without any further restriction on \(S\), one possibility is to compute \(y = \proj(w \mid x + \Delta \B_\infty)\), \(z_i\) for all \(i = 1, \ldots, n\) and retain the \(k\) largest entries, as those will yield the largest sum.
  % \begin{enumerate}
  %   \item form \(w^2\), the vector whose components are \(w_i^2\), \(i = 1, \ldots, n\);
  %   \item compute \(y := \proj(w \mid x + \Delta \B_\infty)\);
  %   \item form \(w - y\) and \((w - y)^2\);
  %   \item form \(z := w^2 - (w - y)^2\);
  %   \item let \(\pi\) be a permutation that sorts the components of \(z\) in decreasing order of value (not magnitude);
  %   \item set \(S_\star\) to \(\{\pi^{-1}(1), \ldots, \pi^{-1}(k)\}\);
  %   \item set \(y_{S_\star^c} = 0\),
  % \end{enumerate}
  % where \(\pi^{-1}\) denotes the inverse permutation.
  Applying the procedure described in \Cref{alg:projB0Binf} with \(L = \varnothing\) corresponds to the steps just outlined.
  By \(\pi^{-1}(1)\), we mean the element of \(F\) that is permuted to first position in the ordering.
  The main cost is the computation of \(\pi\), which can be obtained in \(O(n \log(n))\) operations.

  \begin{algorithm}[ht]
    \caption{%
      \label{alg:projB0Binf}
      Compute the projection of \(w\) into \(C := k \B_0 \cap (x + \Delta \B_\infty)\).
    }
    \begin{algorithmic}[1]
      \Require \(w \in \R^n\), \(w \not \in x + \Delta \B_\infty\), \(L \subseteq \{1, \ldots, n\}\), \(|L| \leq k\)
      \Comment{\(\supp(\proj(w \mid C))\) must contain \(L\)}
      \State compute \(y := \proj(w \mid x + \Delta \B_\infty)\)
      \If{\(|L| = k\)}
        \Return \(L\) and \(\proj(y \mid A_L)\) \Comment{\Cref{lem:suppx,lem:proj-suppx}}
      \EndIf
      \State set \(F := L^c\) and form \(w_F^2\), \(w_F - y_F\), \({(w_F - y_F)}^2\), and \(z := w_F^2 - {(w_F - y_F)}^2\) \Comment{componentwise}
      \State compute a permutation \(\pi\) that sorts the components of \(z\) in decreasing order
      \State set \(S := L \cup \{\pi^{-1}(1), \ldots, \pi^{-1}(k - |L|)\}\) \Comment{\(L\) and the indices of the \(k - |L|\) largest elements of \(z\)}
      \State set \(y_{S^c} = 0\)
      \State \Return \(S\) and \(y\).
    \end{algorithmic}
  \end{algorithm}

  Consider now the case where \(\ell(x) \neq \varnothing\).
  If \(i \in \ell(x)\), \(x + \Delta B_\infty\) cannot intersect any \(A_S\) such that \(i \not \in S\).
  Indeed, any \(y \in \R^n\) such that \(y_i = 0\) satisfies \(\|y - x\|_\infty \geq |y_i - x_i| = |x_i| > \Delta\).
  If \(s(x) = \varnothing\), we are in the context of \Cref{lem:proj-suppx}.
  Thus, we may focus on the case where both \(s(x)\) and \(\ell(x)\) are nonempty.
  Necessarily, \(1 < |s(x)| + |\ell(x)| \leq k\) and \(|\ell(x)| < k\).
  Constraining \(S \subseteq \{1, \ldots, n\}\) to contain \(\ell(x)\) leaves \(k - |\ell(x)|\) indices to be chosen among the remaining \(n - |\ell(x)|\), for a total of
  \[
    {n - |\ell(x)| \choose k - |\ell(x)|} =
    \frac{(n - |\ell(x)|)!}{(k - |\ell(x)|)! \, (n - k)!}
  \]
  possibilities.
  Again, it appears as though the complexity of identifying \(S\) is exponential in \(n\) in the worst case.
  However, the only difference with~\eqref{eq:proj-S} is that \(S_\star\) is now constrained to contain \(\ell(x)\).
  It follows that we may apply \Cref{alg:projB0Binf} with \(L = \ell(x)\).
  If \(m := n - |\ell(x)|\), the procedure has \(O(m \log(m)) = O(n \log(n))\) complexity.

  \section{Implementation and Numerical Results}

  We implemented \Cref{alg:projB0Binf} in the Julia language \citep{bezanson-edelman-karpinski-shah-2017} version \(1.7\) as part of the ShiftedProximalOperators package of \citet{baraldi-orban-shifted-proximal-operators-2022}, whose main objective, as the name implies, is to collect proximal operators of nonsmooth terms with one or two shifts, i.e., \(h(x_k + s_j + t)\), with and without a trust-region constraint, where \(x_k\) and \(s_j\) are fixed iterates set during an outer and an inner iteration.
  ShiftedProximalOperators is used inside the RegularizedOptimization package of \citet{baraldi-orban-regularized-optimization-2022}, which implements, among others, the trust-region methods for nonsmooth regularized problems of \citet{aravkin-baraldi-orban-2021,aravkin-baraldi-orban-2021b}.

  We employ \Cref{alg:projB0Binf} to solve~\eqref{eq:proj} inside two trust-region methods for nonsmooth regularized problems of the form~\eqref{eq:minf+h}.
  The trust region is defined in the \(\ell_\infty\)-norm in both, and provides the box \(x + \Delta \B_\infty\), where \(x\) is the current iterate and \(\Delta\) the trust-region radius.
  At iteration \(j\) of the method of \citet{aravkin-baraldi-orban-2021}, a step \(s_j\) is computed as an approximate solution of the model
  \[
    \minimize{s} \ q(s) + \psi(s; x_j) + \chi(s \mid \Delta \B_\infty),
    \qquad
    q(s) := \nabla f(x_j)^T s + \tfrac{1}{2} s^T B_j s,
  \]
  where \(B_j = B_j^T \in \R^{n \times n}\) is a limited-memory LBFGS or LSR1 approximation of the Hessian of \(f\), and \(\psi(s; x_j) \approx h(x_j + s)\).
  Below, we choose \(\psi(s; x_j) = h(x_j + s) = \chi(x_j + s \mid k \B_0)\) for an appropriate value of \(k \in \N\).
  \(s_j\) is computed using an adaptive stepsize variant of the proximal gradient algorithm named R2 \citep{aravkin-baraldi-orban-2021} that generates inner iterates \(s_{j,l}\), starting with \(s_{j,0} := s_j\).
  At iteration \(l\) of R2, we compute a step \(t_l\) that solves
  \[
    \minimize{t} \ \nabla q(s_{j,l-1})^T t + \tfrac{1}{2} \sigma_l \|t\|_2^2 + \psi(s_{j,l-1} + t; x_j) + \chi(s_{j,l-1} + t \mid \Delta \B_\infty),
  \]
  where \(\sigma_l > 0\).
  If we complete the square and perform the change of variables \(y = x_j + s_{j,l-1} + t\), we obtain a problem of the form~\eqref{eq:proj}.
  We refer to the method outlined above as TR.

  The second trust-region method is a variant specialized to the case \(f(x) = \tfrac{1}{2} \|F(x)\|_2^2\), where \(F: \R^n \to \R^m\) inspired from the method of \citet{levenberg-1944} and \citet{marquardt-1963}, where we redefine \(q(s) := \tfrac{1}{2} \|J(x) s + F(x)\|_2^2\).
  We refer to the latter as LMTR.

  In both methods, the decrease in the model achieved by \(s_j\) is denoted \(\xi\).
  Of particular interest is the decrease achieved by \(s_{j,1}\)---the first step in the inner iterations---which is denoted \(\xi_1\).
  It is possible to show that \(\sqrt{\xi_1}\) may be used as a criticality measure for~\eqref{eq:minf+h}.
  Each method stops as soon as \(\sqrt{\xi_1} \leq \epsilon + \epsilon \sqrt{\xi_{1,0}}\) where \(\xi_{1,0}\) is the \(\xi_1\) observed at the first outer iteration and \(\epsilon = 10^{-6}\).

  We illustrate the behavior of the trust-region methods on the LASSO / basis pursuit denoise problem, in which we fit a linear model to noisy observations \(Ax \approx b\), where the rows of \(A \in \R^{m \times n}\) are orthonormal.
  We set \(b = A x_\star + \varepsilon\), where \(\|x_\star\|_0 = k\) with its nonzero components set to \(\pm 1\) randomly and \(\varepsilon \sim \mathcal{N}(0, 0.01)\).
  In our experiment, we set \(m = 200\), \(n = 512\), and \(k = 10\).
  We formulate the problem as
  \begin{equation}%
    \label{eq:bpdn}
    \minimize{x \in \R^n} \ \tfrac{1}{2} \|Ax - b\|_2^2 + \chi(x \mid k \B_0).
  \end{equation}
  We report results in the form of the solver output in \Cref{lst:BPDN-TR-LSR1,lst:BPDN-TR-LBFGS,lst:BPDN-LMTR}, where \emph{outer} is the outer iteration counter \(j\), \emph{inner} is the number of inner R2 iterations at each outer iteration, \(f(x)\) and \(h(x)\) are the value of the smooth and nonsmooth part of the objective, respectively, \(\sqrt(\xi_1)\) is our criticality measure, \(sqrt{\xi}\) is the square root of the decrease achieved the by step \(s_j\), \(\rho\) is the ratio of actual versus predicted reduction used to accept or reject the step, \(\Delta\) is the trust-region radius, \(\|x\|\) and \(\|s\|\) are the \(\ell_\infty\)-norm of the iterate and step, respectively, \(\|B_j\|\) is the spectral norm of \(B_j\), and \(1 / \nu\) is the regularization parameter \(\sigma_l\) in the R2 model.
  In \Cref{lst:BPDN-TR-LSR1}, \(B_j\) is a limited-memory SR1 operator with memory \(5\).
  In \Cref{lst:BPDN-TR-LBFGS}, \(B_j\) is a limited-memory BFGS operator with memory \(5\).
  All methods use the initial guess \(x_0 = 0\).

  \begin{jllisting}[caption=\label{lst:BPDN-TR-LSR1} TR iterations with L-SR1 on~\eqref{eq:bpdn}.]
outer    inner     f(x)     h(x)     √ξ1      √ξ        ρ       Δ     ‖x‖     ‖s‖    ‖Bⱼ‖
    1        2  1.9e+00  0.0e+00 8.9e-01 8.9e-01  1.5e+00 1.0e+00 0.0e+00 4.7e-01 1.0e+00
    2        9  7.4e-01  0.0e+00 4.8e-01 7.2e-01  1.2e+00 1.4e+00 4.7e-01 5.4e-01 1.0e+00
    3       12  1.0e-01  0.0e+00 1.8e-01 2.3e-01  1.4e+00 1.6e+00 1.0e+00 3.3e-01 1.0e+00
    4       17  3.0e-02  0.0e+00 8.7e-02 1.4e-01  1.0e+00 1.6e+00 1.1e+00 2.9e-01 1.0e+00
    5       22  1.0e-02  0.0e+00 1.6e-02 2.6e-02  1.0e+00 1.6e+00 1.0e+00 3.3e-02 1.0e+00
    6       18  9.5e-03  0.0e+00 2.7e-03 4.4e-03  1.0e+00 1.6e+00 1.0e+00 6.0e-03 1.0e+00
    7        8  9.4e-03  0.0e+00 3.6e-04 3.7e-04  1.5e+00 1.6e+00 1.0e+00 2.6e-04 1.0e+00
    8       10  9.4e-03  0.0e+00 2.0e-04 3.0e-04  1.0e+00 1.6e+00 1.0e+00 3.5e-04 1.0e+00
    9        6  9.4e-03  0.0e+00 2.0e-05 3.3e-05  1.0e+00 1.6e+00 1.0e+00 4.5e-05 1.0e+00
   10        1  9.4e-03  0.0e+00 2.1e-06 2.4e-06  1.2e+00 1.6e+00 1.0e+00 2.7e-06 1.0e+00
TR: terminating with ξ1 = 1.3038641262246793e-6
TR relative error
   norm(TR_out.solution - sol) / norm(sol) = 0.014710272483962346
  \end{jllisting}

  \begin{jllisting}[caption=\label{lst:BPDN-TR-LBFGS} TR iterations with L-BFGS on~\eqref{eq:bpdn}.]
outer    inner     f(x)     h(x)     √ξ1      √ξ        ρ       Δ     ‖x‖     ‖s‖    ‖Bⱼ‖
    1        2  1.9e+00  0.0e+00 8.9e-01 8.9e-01  1.5e+00 1.0e+00 0.0e+00 4.7e-01 1.0e+00
    2       18  7.4e-01  0.0e+00 3.6e-01 7.1e-01  1.0e+00 1.4e+00 4.7e-01 5.4e-01 1.7e+00
    3       23  2.1e-01  0.0e+00 7.1e-02 1.0e-01  1.6e+00 1.6e+00 1.0e+00 9.0e-02 1.8e+00
    4       14  2.0e-01  0.0e+00 4.3e-02 2.7e-01  1.6e+00 1.6e+00 1.0e+00 3.6e-01 2.1e+00
    5       12  8.6e-02  0.0e+00 1.1e-01 2.4e-01  1.2e+00 1.6e+00 1.1e+00 5.0e-01 2.3e+00
    6       18  1.6e-02  0.0e+00 2.9e-02 6.6e-02  1.2e+00 1.6e+00 1.1e+00 1.3e-01 2.5e+00
    7       23  1.1e-02  0.0e+00 1.4e-02 2.7e-02  1.4e+00 1.6e+00 1.1e+00 3.2e-02 2.6e+00
    8       20  9.8e-03  0.0e+00 7.4e-03 1.7e-02  1.1e+00 1.6e+00 1.0e+00 2.4e-02 2.6e+00
    9       14  9.5e-03  0.0e+00 1.7e-03 3.7e-03  1.2e+00 1.6e+00 1.0e+00 5.3e-03 2.7e+00
   10       14  9.4e-03  0.0e+00 7.7e-04 1.6e-03  1.2e+00 1.6e+00 1.0e+00 1.6e-03 2.5e+00
   11       15  9.4e-03  0.0e+00 3.0e-04 6.1e-04  1.2e+00 1.6e+00 1.0e+00 6.1e-04 2.4e+00
   12        9  9.4e-03  0.0e+00 8.9e-05 2.0e-04  1.2e+00 1.6e+00 1.0e+00 3.1e-04 2.5e+00
   13        8  9.4e-03  0.0e+00 3.1e-05 6.2e-05  1.3e+00 1.6e+00 1.0e+00 7.6e-05 2.5e+00
   14        8  9.4e-03  0.0e+00 1.4e-05 3.0e-05  1.2e+00 1.6e+00 1.0e+00 3.0e-05 2.5e+00
   15        4  9.4e-03  0.0e+00 4.3e-06 8.6e-06  1.1e+00 1.6e+00 1.0e+00 5.5e-06 2.6e+00
   16        3  9.4e-03  0.0e+00 2.5e-06 5.8e-06  1.0e+00 1.6e+00 1.0e+00 4.3e-06 2.6e+00
TR: terminating with ξ1 = 1.0999297328606739e-6
TR relative error
   norm(TR_out.solution - sol) / norm(sol) = 0.014709629662551134
  \end{jllisting}

  \begin{minipage}{\linewidth}  % to prevent a page break between the caption and listing.
    \begin{jllisting}[caption=\label{lst:BPDN-LMTR} LMTR iterations on~\eqref{eq:bpdn}.]
outer    inner     f(x)     h(x)     √ξ1      √ξ        ρ       Δ     ‖x‖     ‖s‖     1/ν
    1        9  1.9e+00  0.0e+00 8.9e-01 1.4e+00  1.0e+00 1.0e+00 0.0e+00 1.0e+00 1.0e+00
    2       11  1.1e-02  0.0e+00 2.3e-02 4.1e-02  1.0e+00 3.0e+00 1.0e+00 7.0e-02 1.0e+00
    3       11  9.4e-03  0.0e+00 3.8e-04 6.8e-04  1.0e+00 3.0e+00 1.0e+00 1.2e-03 1.0e+00
    4        4  9.4e-03  0.0e+00 7.0e-06 1.2e-05  1.0e+00 3.0e+00 1.0e+00 1.8e-05 1.0e+00
LMTR: terminating with ξ1 = 2.797637121965124e-12
LMTR relative error
     norm(LMTR_out.solution - sol) / norm(sol) = 0.014710437655962767
    \end{jllisting}
  \end{minipage}

\begin{figure}[ht]
    \includetikzgraphics[width=.32\linewidth]{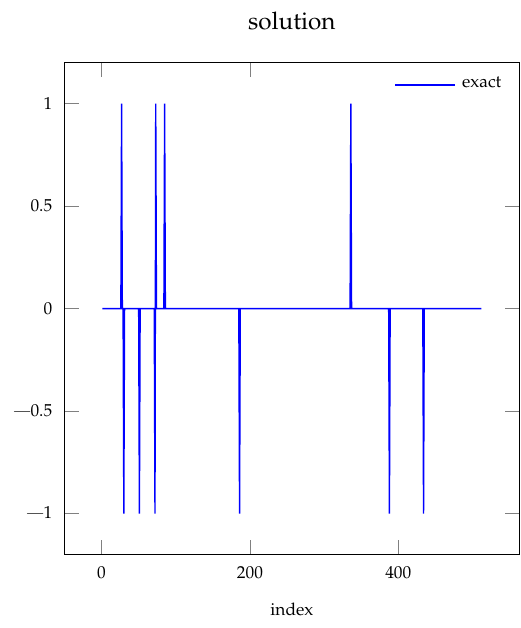}
    \hfill
    \includetikzgraphics[width=.32\linewidth]{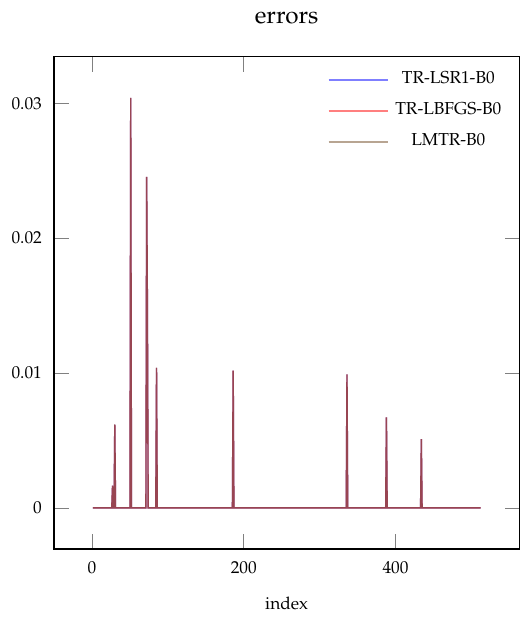}
    \hfill
    \includetikzgraphics[width=.32\linewidth]{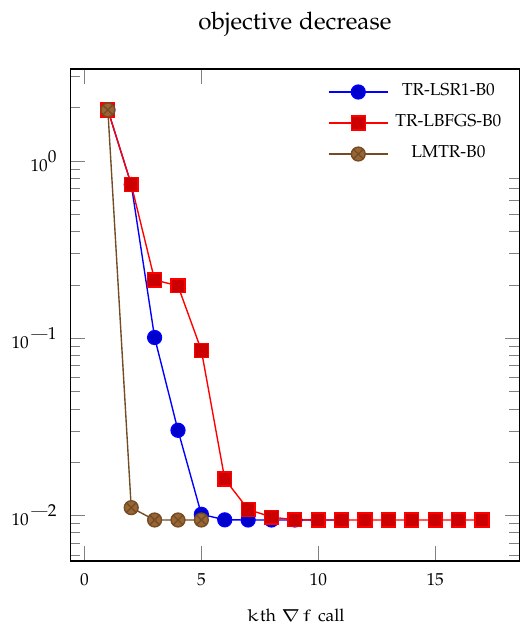}
    \caption{%
      \label{fig:bpdn}
      Exact solution of~\eqref{eq:bpdn} (left), absolute errors (center), and objective decrease history as a function of the number of \(\nabla f\) evaluations (right).
    }
  \end{figure}

  \Cref{fig:bpdn} shows the exact solution \(x_\star\), and the objective history of each solver.
  All three solvers find a solution where the amplitude of the peaks are within \(10^{-2}\) of the correct amplitude.
  It is not surprising that LMTR, which exploits the least-squares structure of~\eqref{eq:bpdn} performs better than TR; its model is exact at each iteration, which is reflected in the fact that \(\rho = 1\) at each iteration in \Cref{lst:BPDN-LMTR}.
  TR also performs well, although, surprisingly, the potentially indefinite L-SR1 Hessian approximations of the positive definite Hessian \(A^T A\) yield fewer iterations than the positive-definite L-BFGS approximation.

  From a computation cost point of view, each outer TR iteration costs one evaluation of \(f\) and, if the step is accepted, one evaluation of \(\nabla f\).
  In \Cref{lst:BPDN-TR-LSR1,lst:BPDN-TR-LBFGS}, every step is accepted.
  Each inner R2 iteration in TR costs a product between the limited-memory quasi-Newton approximation and a vector, and an execution of \Cref{alg:projB0Binf}.
  Each outer LMTR iteration costs one evaluation of \(F(x)\).
  Each inner R2 iteration in LMTR costs a Jacobian-vector product, a transposed-Jacobian-vector product, and an executation of \Cref{alg:projB0Binf}.
  % Of course, we could have solved~\eqref{eq:bpdn} directly with R2 without recourse to \Cref{alg:projB0Binf}, but that requires~\(35\) iterations, i.e., \(35\) gradient evaluations.

\begin{figure}[ht]
    \begin{center}
      \includetikzgraphics[width=.32\linewidth]{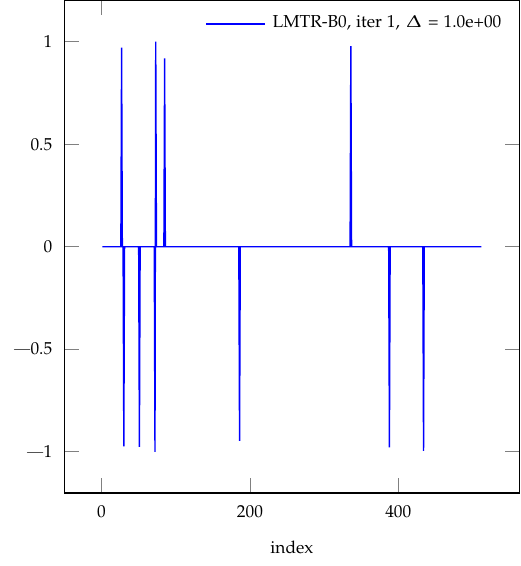}
      \hfill
      \includetikzgraphics[width=.32\linewidth]{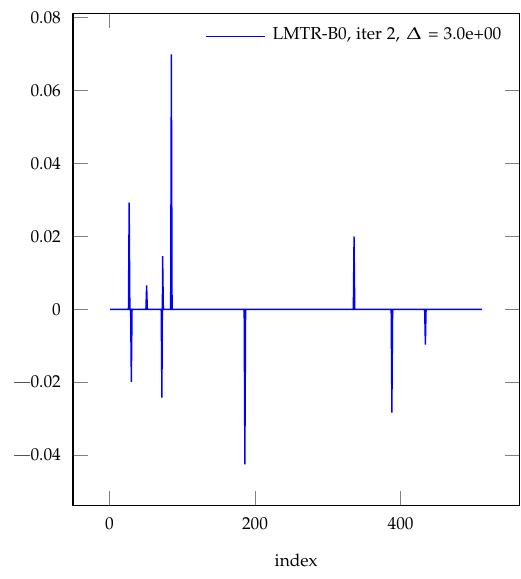}
      \hfill
      \includetikzgraphics[width=.32\linewidth]{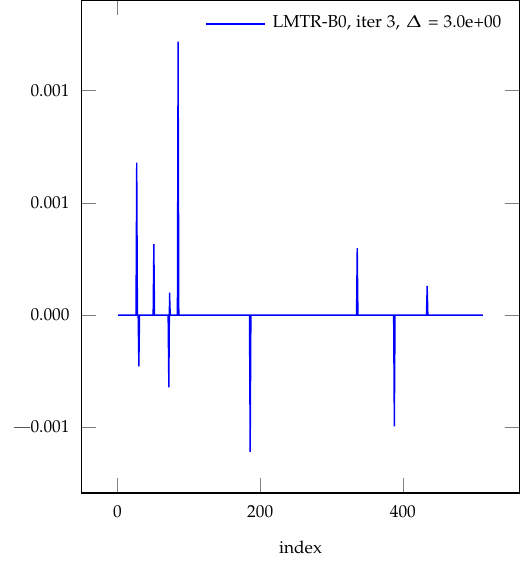}
      % \hfill
      % \includetikzgraphics[width=.2\linewidth]{bpdn-steps-LMTR-B0-4}
      % \hfill
      % \includetikzgraphics[width=.2\linewidth]{bpdn-steps-LMTR-B0-5}
    \end{center}
    \caption{%
      \label{fig:lmtr-steps}
      First three steps generated during the iterations of LMTR applied to~\eqref{eq:bpdn}.
      At iteration~\(1\), the trust-region constraint is active (left).
      It is inactive at subsequent iterations.
    }
  \end{figure}

  In each method, each step is a sum of R2 steps, each of which is a projection of the form~\eqref{eq:proj}.
  \Cref{fig:lmtr-steps} shows the first three LMTR steps.
  At iteration~\(1\) (leftmost plot), the trust-region constraint is active, i.e., the step norm \(\|s\|_\infty = \Delta\), which means that at least one of the projections computed during the R2 iterations resulted in a point in \(k \B_0\) at the boundary of \(x + \Delta \B_\infty\).
  At subsequent LMTR iterations, \(\|s\|_\infty < \Delta\), which is expected in trust-region methods as convergence occurs, and means that at least the final projection computed during the R2 iterations resulted in a point lying strictly inside \(x + \Delta \B_\infty\).

  \section{Closing Remarks}

  Although \(C\) is a nonconvex set, there exists an efficient projection into it, and the latter can be used to design proximal methods for nonsmooth regularized problems \citep{aravkin-baraldi-orban-2021,aravkin-baraldi-orban-2021b}.
  \Cref{alg:projB0Binf} makes it possible to solve sparsity-constrained problems by way of trust-region methods.
  It also makes it conceivable to tackle the more general problem~\eqref{eq:min-sparse-box} by way of one of the algorithms proposed by \citep{beck-eldar-2013}.

  Possible extensions of this work include balls defined by other norms, such as other \(\ell_p\) norms or elliptical norms.
  However, it is not clear that \Cref{alg:projB0Binf} generalizes in a straightforward way.
  Indeed, the key is that the projection into \(x + \Delta \B_\infty\) is defined componentwise.
  It is not difficult to sketch an example where the same procedure using the Euclidean norm yields an erroneous projection.

  Another possible generalization is to consider \(x \not \in k\B_0\), as might occur in an infeasible method.

  The exploration of such generalizations is the subject of ongoing research.

  \subsection*{Acknowledgements}

  The author wishes to thank Aleksandr Aravkin and Robert Baraldi for fruitful discussions that made this research possible. % and improved the presentation.

  \small
  \bibliographystyle{abbrvnat}
  \bibliography{abbrv,b0-binf-prox}

\end{document}